\documentclass{amsart}
\usepackage{graphicx, amsfonts, amssymb}
\usepackage{mathrsfs, amsmath, amsthm}
\usepackage{verbatim}
\usepackage{relsize}
\usepackage{epsfig}
\usepackage{color, colortbl}
\usepackage{comment}
\usepackage{mymacros}
\usepackage{xfrac,bigints}

\usepackage{hyperref}
\hypersetup{
	colorlinks=true,
	linkcolor=blue,
	filecolor=magenta,      
	urlcolor=cyan,
	pdftitle={Generalized Hilbert operators arising from Hausdorff matrices},
	pdfpagemode=FullScreen,
}

\newtheorem{theorem}{Theorem}[section]

\theoremstyle{definition}

\theoremstyle{remark}
\newtheorem{remark}[theorem]{Remark}

\numberwithin{equation}{section}

\def\a{\alpha}                  \def\g{\gamma}
            
\def\z{\zeta}                  \def\th{\theta}

\def\m{\mu}

\def\G{\Gamma}

\def\T{{\mathbb T}}
  \def\N{{\mathbb N}}

\def\H{{\mathcal H}}
\newcommand{\C}{\mathbb C}
\newcommand{\D}{\mathbb D}

\def\cC{{\mathcal C}}    \def\cd{{\mathcal D}}
\def\cK{{\mathcal K}}

\newcommand{\n}[1]{\Vert#1\Vert}

\newcommand{\abs}[1]{\lvert#1\rvert}

\title{Generalized Hilbert operators arising from Hausdorff matrices}

\author{C. Bellavita }
\email{cbellavita@math.auth.gr}
\address{Department of Mathematics, Aristotle University of Thessaloniki, 54124, Thessaloniki, Greece.}

\author{N. Chalmoukis}
\email{nikolaos.chalmoukis@unimib.it}
\address{Dipartimento di Matematica e Applicazioni, Universit\'a degli studi di Milano Bicocca, via Roberto Cozzi, 55 20125, Milano, Italy}

\author{V. Daskalogiannis}
\email{vdaskalo@math.auth.gr}
\address{Department of Mathematics, Aristotle University of Thessaloniki, 54124, Thessaloniki, Greece.}

\author{G. Stylogiannis}
\email{stylog@math.auth.gr}
\address{Department of Mathematics, Aristotle University of Thessaloniki, 54124, Thessaloniki, Greece.}

\thanks{This research project was supported by the
Hellenic Foundation for Research and Innovation (H.F.R.I.) under the '2nd
Call for H.F.R.I. Research Projects to support Faculty Members \& Researchers' (Project Number: 4662).}



\begin{document}
\begin{abstract}
For a finite, positive, Borel measure $\mu$ on $(0,1)$ we consider an infinite matrix $\G_\mu$, related to the classical Hausdorff matrix defined by the same measure $\m$, in the same algebraic way that the  Hilbert  matrix is related to the Ces\'aro matrix. When $\mu$ is the Lebesgue measure, $\G_\mu$ reduces to the classical Hilbert matrix. We prove that the matrices $\G_\mu$ are not Hankel, unless $\mu$ is a constant multiple of the Lebesgue measure,
 we give necessary and sufficient conditions for their boundedness on the scale of Hardy spaces $H^p, \, 1 \leq p <  \infty$, and we study their compactness  and complete continuity properties. In the case $2\leq p<\infty$, we are able to compute the exact value of the norm of the operator.
\end{abstract}

\maketitle
\section{Introduction}

Let $\D$ be the unit disc in the complex plane $\C$, and let $H(\D)$ be the Fr\'echet space of all analytic functions on $\D$.
Let   $0<p<\infty$ and   $f\in H(\mathbb{D})$. For  $0\leq r<1$ let
\[
M_p(r,f)\,:=\,\left(\dfrac{1}{2\pi}\int_{0}^{2\pi}\vert f(re^{i\theta})\vert^p\, d\theta\right)^{\frac{1}{p}}\,
\]
be the usual $p$-integral means of $f$ on $|z|=r$.  The
Hardy space $H^p =H^p(\mathbb{D})$ consists of the functions $f\in H(\D)$ for which
\[
\norm{f}_{H^p}\,:=\,\sup_{0\leq r<1}M_p(r,f) <+\infty,
\]
while for $p=\infty$, $\,H^\infty(\D)$ consists of bounded analytic functions on $\mathbb{D}$, i.e.
\[
\n{f}_{\infty}:=\sup_{z\in\D}|f(z)|<\infty.
\]
For $1\leq p \leq \infty$ the above quantities are norms and $H^p$ are  Banach spaces.
We will use some basic properties of these spaces
which we state briefly. If $0 < p < q< \infty$ then $ H^p\supset H^q\supset H^{\infty}$.
If $f\in H^p$, the  growth estimate \cite[p. 36]{Duren1970} holds,
\begin{equation}\label{growth}
\vert f(z)\vert
\,\leq\,
 \left(\dfrac{2}{1-\vert z\vert}\right)^{\frac{1}{p}}\norm{f} _{H^p},\quad z\in \mathbb{D}.
\end{equation}
Each $f\in H^p$ has radial limits $f^{*}(\z)=\lim_{r\to 1^{-}}f(r\z)$ for a.e. $\z\in \bT:=\partial\D$ with respect to
the  Lebesgue measure $|d\z|= d\theta, \,\, \z = e^{i\theta}$. 
The boundary function $f^{*}$ is $p$-integrable with $\norm{f^{*}}_{L^p(\bT)}=\norm{f}_{H^p}$.
In the sequel we use $f$, instead of
$f^{*}$, to denote the boundary function.
For $1<p<\infty$ the dual space $(H^p)^*$ can be identified isomorphically with $H^q,\;\frac{1}{p}+\frac{1}{q}=1$,
under the Cauchy pairing
\[
\langle f, g \rangle \,=\, \frac{1}{2\pi}\int_\bT f(\z)\,\overline{g(\z)}\,|d\z| ,\quad f\in H^p,\;g\in H^q.
\]

\subsection*{The Ces\'aro and Hilbert matrices.}
Consider the Ces\'aro matrix $C$  and the Hilbert matrix $H$,
\[
C=\left(
\begin{array}{ccccc}
            1 & 0  & 0  & .  \\ [4pt]
  \frac{1}{2} & \frac{1}{2}  & 0  & .  \\ [4pt]
  \frac{1}{3} & \frac{1}{3}  & \frac{1}{3}  & . \\
  .  & . & . & .  \\
\end{array}%
\right), \quad \quad H=\left(
\begin{array}{ccccc}
            1 & \frac{1}{2}  & \frac{1}{3}  & .  \\ [4pt]
  \frac{1}{2} & \frac{1}{3}  & \frac{1}{4}  & .  \\ [4pt]
  \frac{1}{3} & \frac{1}{4}  & \frac{1}{5}  & . \\
  .  & . & . & .  \\
\end{array}
\right).
\]
These matrices represent two of the most well-studied bounded linear operators on the Hardy space $H^2$, which we will denote by $\cC$ and $\cH$ respectively, and they 
act  on  functions $f(z)=\sum_{n=0}^{\infty}a_nz^n\in H^2$ as follows;
\begin{equation}\label{Ces-m}
\cC f(z)\,:=\,\sum_{n=0}^{\infty}\left(\frac{1}{n+1}\sum_{k=0}^na_k\right)z^n\,,
\end{equation}
and
\begin{equation}\label{Hil-m}
\cH f(z)
:=\sum_{n=0}^{\infty}\left(\sum_{k=0}^{\infty}\frac{a_k}{n+k+1}\right)z^n.
\end{equation}
Note that the series for $\cC f$ is well defined and converges on $\D$ for all $f\in H(\D)$,  thus $\cC f\in H(\D)$,  while the series for $\cH f$ is not always defined. For example, when  $f(z)=\frac{1}{1-z}$ the coefficients in the formal power series \eqref{Hil-m} are not defined. However,  Hardy's inequality \cite[p. 48]{Duren1970}, guarantees that
\[
\sum_{n\geq 0}\dfrac{|a_n|}{n+1}<\infty,
\]
for  functions in the Hardy space $H^1$. In particular, this implies  that $\cH f$ is a well defined  analytic function on $\D$ for each $f\in H^1$.

A calculation shows that $\cC f$ and $ \cH f$ can be written in integral form,
\begin{equation}\label{Ces-g}
\begin{split}
 \cC f(z)\,&=\,\frac{1}{z}\int_0^zf(\zeta)\frac{1}{1-\z}\,d\zeta=\int_0^1f(tz)\frac{1}{1-tz}\,dt
 \\
 &=\int_0^1f(tz)g'(tz)\,dt,
 \end{split}
\end{equation}
and
\begin{equation}\label{Hil-g}
\cH f(z)=\int_0^1f(t)\frac{1}{1-tz}\,dt=\int_0^1f(t)g'(tz)\,dt,
\end{equation}
where in both cases $g(z)=\log\dfrac{1}{1-z}$. The convergence of the integral for  $\cH f$ is  guaranteed
for all $f\in H^1$ by the Fej\'er-Riesz inequality \cite[Theorem 3.13]{Duren1970}.
It is well known that $\cC$ is a bounded linear operator on all  $H^p,\;0< p<\infty$, \cite{Siskakis1987}, \cite{Miao1992} and that
$\cH$ is bounded on $H^p$ for $1< p<\infty$, but not bounded on $H^1$ and $H^\infty$ \cite{Diamantopoulos2000}.

Observe that  $H$ is obtained from $C$ by the following formal manipulation;
\textit{Eliminate  the zeros in each column of  $C$ and
shift up the columns to their  first nonzero entry.}
In rigorous terms this is equivalent to the following algebraic relation. Let $e_n(z)=z^n,  n=0,1,2, ...$ be the monomials, which 
form an orthonormal basis of $H^2$. Then, we have
 $$
\cC e_n(z)=e_n(z)\cH e_n(z) = \frac{1}{z}\left(\log\frac{1}{1-z}- \sum_{k=1}^n\frac{z^k}{k}\right).
 $$
Equivalently,
\begin{equation}\label{H from C}
\H e_n(z)=\frac{1}{z^n} \cC e_n(z).
\end{equation}
Note that this can be expressed in terms of the shift operator $Sf(z):=zf(z)$ as
$$
\cC e_n=S^n\H e_n, \,\,  n=0,1,2,...
$$

\subsection*{Generalized Hilbert matrices arising from Hausdorff matrices}\label{Haus}
For a given sequence $\mu=\{\mu_n\}_{n=0}^\infty$ the Hausdorff matrix induced by $\mu$ is the
 lower triangular matrix
 \begin{equation}
\cK_{{\m}}= \left(
\begin{matrix}
c_{00}& 0 & 0&  \cdots     \\
c_{10}& c_{11} &0& \cdots     \\
c_{20} & c_{21} & c_{22} & \cdots    \\
\vdots &\vdots  &\vdots     &\ddots
\end{matrix}
\right)
\notag
\end{equation}
with entries $c_{nk}$ given by
\begin{equation}
c_{n k}=\binom{n}{k}\Delta^{n-k}\mu_k, \quad     0\leq k\leq n,
\notag
\end{equation}
where $\Delta\mu_n=\mu_n-\mu_{n+1}$ is  the forward
difference operator, with iterates
$$
\Delta^{k}\mu_n=\Delta(\Delta^{k-1}\mu_n),
\quad  k=1,2,\dots, \quad \Delta^0\mu_n=\mu_n.
$$
 In the special case where $\{\mu_n\}$ is the moment sequence
\begin{equation}
\mu_n=\int_0^1t^n\,d\mu(t), \qquad n=0,1,\dots
\notag
\end{equation}
of a finite positive Borel measure $d\mu$ on $(0,1)$, the entries of $\cK_{\m}$ can be calculated explicitely; 
\begin{equation}
c_{n k}=\binom{n}{k}\int_0^1t^k(1-t)^{n-k}\,d\mu(t),
\qquad 0\leq k\leq n\,.
\notag
\end{equation}
Hausdorff matrices have been studied in connection to the classical summability
theory on sequence  spaces \cite{Hardy1943} and recently on spaces of analytic functions
 \cite{Galanopoulos2001}, \cite{Galanopoulos2006}.

The Ces\'aro matrix is a typical example in the family of  Hausdorff  matrices, and is obtained by the
moment sequence $\mu_n=\frac{1}{n+1}$ of the Lebesgue measure, $d\mu(t)=dt$. Applying the
operation of erasing zeros and shifting up the columns, i.e. the analogous of (\ref{H from C}) to the matrix $\cK_{\mu}$, we obtain a new matrix
$$
\G_{\m}=\left(%
\begin{array}{ccccc}
            \g_{00} & \g_{01}  & \g_{02}  & .  \\
\g_{10} & \g_{11}  & \g_{12}  & .  \\
  \g_{20} & \g_{21}  & \g_{22}  & . \\
  .  & . & . & .  \\
\end{array}%
\right)
$$
with entries
$$
\g_{n k}=c_{n+k, k}=\binom{n+k}{k}\int_0^1t^k(1-t)^n\,d\m(t)\,.
$$
By its  construction $\G_{\mu}$ reduces to the Hilbert matrix $H$ for the
particular choice of the measure  $d\m(t)=dt$, i.e.  this choice of $\mu$  gives
$\g_{n,k}=\frac{1}{n+k+1}$. In a very recent article \cite{Athanasiou2023}, the author studied the boundedness properties  of the matrix $\G_{\mu}$, when acting on the classical sequence spaces $\ell^p, 1\leq p \leq \infty$. However, a different approach is required when one studies the operator on spaces of holomorphic functions.

In this article we will study the action of the operators $\G_{\m}$ on the Hardy spaces $H^p,\, 1\leq p < \infty$.
If  $f(z)=\sum_{k=0}^{\infty} a_kz^k$ is analytic on the unit disc then we define the formal power series
\[
\G_{\m} f(z)=\sum_{n=0}^{\infty}A_nz^n,\quad\text{where} \,\,\, A_n=\sum_{k=0}^{\infty}\g_{n,k}a_k, \quad n=0,1, ...,
\]
whenever the series defining the coefficients $A_n$ are convergent.  If $f$ is a polynomial then the sum giving
$A_n$  is finite for each $n$ and one can verify that the formal power series for $\Gamma_\mu f$ converges in the unit disc.  Therefore, if one can establish an inequality of the type 
\[ \norm{\Gamma_\mu f}_{H^p} \leq C_p \norm{f}_{H^p}\,, \]
for some $C_p>0$ and for all polynomials $f$, then, by the density of the polynomials in $H^p$, the operator $\Gamma_\mu$ extends in a unique way to a bounded linear operator on the whole $H^p$.


\section{An integral representation and main results}

The integral giving $\H$ in (\ref{Hil-g})  can be viewed as an \lq\lq improper\rq\rq\, line integral, with the path of
integration being the radius $[0, 1]$. When $\H$ is applied to a function belonging to $H^1$, it turns
out that we can change the path of integration in \eqref{Hil-g} to be the arc
$$
\gamma(s)= \gamma_z(s)= \frac{s}{1-(1-s)z}, \quad 0\leq s\leq 1\,,
$$
which joins $0$ to $1$ and lies inside $\D$ for $0\leq s<1$ for every $z\in \D$,
see \cite[Section 2]{Diamantopoulos2000}. With this change of variable we obtain
\begin{align*}
\H f(z)=\int_0^1\frac{1}{1-(1-t)z}f\left(\frac{t}{1-(1-t)z}\right)\,dt
=:\int_0^1T_t f(z)\,dt\,.
\end{align*}
This gives a representation of $\cH$ in terms of the weighted composition operators

\begin{equation}\label{T-t}
T_t(f)(z)=w_t(z)f(\phi_t(z))\,dt\,,
\end{equation}
with
\begin{equation}\label{w-phi}
\phi_t(z):=\frac{t}{1-(1-t)z}, \quad             w_t(z):= \frac{1}{1-(1-t)z}=\frac{\phi_t(z)}{t}\,.
\end{equation}
It is easy to check that each $T_t$, $0<t<1$,  is  bounded on $H^p$ for $0< p<\infty$.

We will derive an analogous representation for $\G_{\m}$. We assume that $\mu$ is a finite positive Borel measure on $(0,1)$. Assume also that $f$ is such that $\G_\mu f$ is well defined as an analytic function in $\mathbb{D}$ (for example when $f$ is a polynomial). In that case, all sums and integrals are interchangeable.
We then have:
\[
\begin{split}
\Gamma_\mu f(z) &= \sum_{n= 0}^{\infty}\left(\sum_{k= 0}^{\infty}a_k\, \binom{n+k}{k}\int_0^1
(1-t)^n t^k\,d\mu(t)\right)\,z^n
\\
&=\int_0^1\,\sum_{k=0}^{\infty} a_k t^k \, \sum_{n=0}^{\infty} \binom{n+k}{k} (1-t)^n z^n\,d\mu(t)\,.
\end{split}
\]
Next, using the identity
\begin{equation}\label{identity}
\sum_{n= 0}^{\infty}\binom{n+m}{n}w^n=\dfrac{1}{(1-w)^{m+1}}\,,
\end{equation}
we get
\[
\begin{split}
\sum_{n=0}^{\infty} \binom{n+k}{n} \left((1-t) z\right)^n &=
\frac{1}{(1-(1-t)z)^{k+1}}\\
&=\left(\frac{1}{(1-(1-t)z)}\right)^k\frac{1}{(1-(1-t)z)}\,.
\end{split}
\]
Therefore,
\[
\begin{split}
\G_\mu f(z) &=\int_0^1\,\sum_{k= 0}^{\infty} a_k t^k\left(\frac{1}{(1-(1-t)z)}\right)^k\frac{1}{(1-(1-t)z)}
\,d\mu(t) \\
&=\int_0^1\sum_{k=0}^{\infty} a_k  \left(\frac{t}{(1-(1-t)z)}\right)^k\frac{1}{(1-(1-t)z)}\,d\mu(t)
\\
&=\int_0^1 f\left(\frac{t}{(1-(1-t)z)}\right)\frac{1}{(1-(1-t)z)}\,d\mu(t)\\
&=\int_0^1w_t(z) f(\phi_t(z))\,d\mu(t)\,,
\end{split}
\]
where $w_t$ and $\phi_t$ are as in (\ref{w-phi}). Hence
\[
\G_\mu f(z)=\int_0^1 T_t f(z) \,d\mu(t)\,,
\]
where $T_t$ are the weighted composition operators (\ref{T-t}).
 This formula can be used to study $\G_{\mu}$ on spaces of analytic functions, exploiting information available for
 (weighted) composition operators. An important property of the operators $T_t$ is the following.

\begin{prop}\label{the adjoint}
Let $T_{t}^{*}$ denote the adjoint of the operator $T_{t}$ on $H^p$, for $1<p<\infty$. Then,
$$
T_{t}^{*}=T_{1-t}\,,
$$
for every $0<t<1$. In particular $T_{1/2}$ is self adjoint in $H^2$.
\end{prop}
\begin{proof} Let $1<p<\infty$ and consider $T_t: H^p\to H^p$. The adjoint $T_t^{*}$ acts on $H^q$ with  $\frac{1}{p}+\frac{1}{q}=1$, and we have that
\[
\int_\T {T_tf(\zeta)}\,\overline{g(\zeta)}\,| d\z |\,=\,\int_{\T} f(\zeta)\,\, \overline{T_t^*g(\zeta)}\,| d\z |\,. 
\]

Assume  that $f(z)=\sum_{k=0}^{\infty}a_k z^k$ and  $g(z)=\sum_{m=0}^{\infty}b_mz^m$ are polynomials. Then,

\[\begin{split}
\int_\T {T_tf(\zeta)}\,\overline{g(\zeta)}\,| d\z |  =&
\int_\T \frac{1}{1-(1-t)\zeta}\, f\left(\frac{t}{1-(1-t)\zeta}\right)\,\overline{g(\zeta)}\,| d\z |
\\
=&\int_{\T} \sum_{k} a_k t^k \, \dfrac{1}{(1-(1-t)\zeta)^{k+1}} \,\, \overline{g(\zeta)}\, | d\z |
\\
=&\int_{\T} \sum_{k} a_k t^k\, \sum_{n } \binom{n+k}{n} \left((1-t) \zeta\right)^n \overline{g(\zeta)} \,| d\z |
\\
=&\int_{\T} \sum_{k} a_k t^k\, \sum_{n} \binom{n+k}{n} \left((1-t) \zeta \right)^n\,\overline{\sum_{m}b_m\z^m}\, | d\z |\,.
\end{split}
\]
Notice here that
\[
\int_\T \z^n \overline{\z}^m \,| d\z |=
\begin{cases}
2\pi, \quad n=m
\\
0,\quad n\neq m
\end{cases}\,,
\]
and also that
\[
\begin{split}
T_{1-t}(g)(z) &= \dfrac{1}{1-tz}\,g\left(\dfrac{1-t}{1-tz} \right)
\\
&=\sum_{m} b_m (1-t)^m\, \dfrac{1}{(1-tz)^{m+1}}
\\
&=\sum_{m} b_m (1-t)^m\, \sum_{n} \binom{n+m}{n}\,t^nz^n\,,
\end{split}
\]
hence
\[
\begin{split}
\int_\T {T_tf(\zeta)}\,\overline{g(\zeta)}\,&| d\z |
=2\pi \sum_{k} a_k \left( \sum_{m} t^k\binom{m+k}{m}\overline{b_m}(1-t)^m\right)
\\
=&\int_{\T} \sum_k a_k \zeta^k \left( \sum_{m}  \binom{m+k}{m}\Bar{b}_m\,(1-t)^m\right) t^k \Bar{\z}^k \,| d\z |
\\
=&\int_{\T} f(\zeta)\, \overline{\sum_n \left(\sum_{m} {b_m}(1-t)^m\,\binom{m+n}{m} \right) t^n\,\zeta^n} \, | d\z |
\\
=&\int_{\T} f(\zeta)\, \overline{ \sum_{m} {b_m}(1-t)^m \,\sum_n  \binom{n+m}{n}\,t^n\,\zeta^n }\, | d\z |
\\
=& \int_{\T} f(\zeta)\,\, \overline{T_{1-t}g(\zeta)}\,| d\z |\, .
\end{split}
\]
The above holds for all polynomials $f\in H^p$ and $g\in H^q$.  Since the polynomials are dense in $H^p$ and $H^q$, we have verified that $T^*_t=T_{1-t}$.
\end{proof}

In \cite{Diamantopoulos2000} the authors estimated the norm of the weighted composition operator $T_t$, and we know that
\[
\n{T_tf}_{H^p}\leq\, \dfrac{1}{t^{1-\frac{1}{p}}\,(1-t)^{\frac{1}{p}}} \;\n{f}_{H^p},\;\;p\geq 2\,.
\]
In Section \ref{proofs} we estimate the norm of $T_t$, considering all $p\geq 1$, providing an alternative proof to the one in \cite{Diamantopoulos2000}, in which we do not have to translate the problem into a problem of the upper half-plane. Our first theorem shows that in the class of the operators $\Gamma_\mu$, being a Hankel operator is the exception rather than the rule. 
\begin{thm}\label{Hankel}
The operator $\Gamma_{\mu}$ is a Hankel operator if and only if  $\mu$ is a constant multiple of the Lebesgue measure.
\end{thm}
 We also note that there is no measure $\mu\neq 0$ for which the operator $\Gamma_{\mu}$ is Toeplitz. To see this, assume that $\Gamma_\mu$ is a Toeplitz operator. Then, for the elements of the main diagonal, we must have that $\gamma_{00}=\gamma_{11}$ or equivalently $1=2\int_{0}^1 t(1-t)d\mu(t)$.
However, $2t(1-t)\leq 1/2$ for each $t\in [0,1]$, which leads to the contradiction $1 \leq \frac{1}{2}$.

 Next we study necessary and sufficient conditions for continuity of the operators $\Gamma_\mu$.
In order to formulate our main result it will be convenient to introduce the following weight function;  
\[ \psi_p(t):= 
\begin{cases}
    \dfrac{t^{\frac{1}{p}-1}}{(1-t)^{\frac{1}{p}}}\,,\,\,\, &\text{if} \,\, p>1
    \\[0.2in]
   \log\left(\dfrac{e}{t}\right)\dfrac{1}{(1-t)}\,, \,\,\,&\text{if}\,\, p=1\,.
\end{cases}
\]
    \begin{thm}\label{main theom p>1}
The operator $\Gamma_\mu$ is bounded on  $H^p, \;1\leq p<\infty$,  if and only if
\begin{equation}\label{condition p>1}
    \int_{0}^{1}\psi_p(t)d\mu(t)<\infty\,.
\end{equation}
Furthermore, there exist positive constants $A_p,\,B_p$ depending only on $p$ such that 
\[
A_p \int_0^1 \psi_p(t)\,d\mu(t)\, \leq \,\norm{\G_\mu}_{H^p \to H^p}\,\leq\, B_{p} \int_{0}^{1} \psi_p(t)\,d\mu(t)\,. 
\]
In particular, when $2\leq p < \infty$, this is the exact value of the norm, i.e.
\[
\|\Gamma_\mu\|_{H^p\to H^p}=\int_{0}^{1}\dfrac{t^{\frac{1}{p}-1}}{(1-t)^{\frac{1}{p}}}d\mu(t)\, .
\]
\end{thm}

Finally, we consider the question of compactness and complete continuity. It turns out that $\G_\mu$ have a behaviour similar to the classical Hilbert matrix, in terms of compactness. Nonetheless, in the endpoint case $p=1$, we prove that the operators $\G_\mu$ are completely continuous. An operator $T:X\to X$ is completely continuous if it maps every relatively weakly compact subset of $X$ into a relatively compact subset of $X$.  In general, every compact operator is completely continuous, however the converse is false when $X$ is non-reflexive.

\begin{thm}\label{Compact}
The operator $\Gamma_{\mu}$ is not compact on $H^p$, $1\leq p<\infty$, unless $\mu $ is the zero measure. However,  $\Gamma_\mu$ is completely continuous on $H^1$ whenever it is bounded.
\end{thm}



\section{Proof of main results}\label{proofs}
We proceed with the proofs of our main results. Without loss of generality, we always assume that $\mu(0,1)=1$.

\begin{proof}[Proof of Theorem \ref{Hankel}]
If $\mu$ is the  Lebesque measure then it is well known that $\Gamma_{\mu}$ is Hankel. For the converse implication, let us suppose that $\Gamma_{\mu}$ is Hankel.
We will use induction to prove that $\mu$ is the Lebesgue measure.
We know that the entries of the matrix $\Gamma_{\mu}$ are
$$
\gamma_{n,k}=\binom{n+k}{n}\int_{0}^{1}t^{n}(1-t)^{k}\,d\mu(t)\,.
$$
For $n=1$ we have
$\gamma_{1,0}=\gamma_{0,1}$, that is
$$
\int_{0}^{1}t\,d\mu=\int_{0}^{1}(1-t)\,d\mu(t)\,.
$$
This implies, since $\mu(0,1)=1$, that
$$
\int_{0}^{1}t\,d\mu(t)=\frac{1}{2}\,.
$$
Let us suppose that
\begin{equation}\label{momentdm}
    \int_{0}^{1}t^{m}\,d\mu(t)=\frac{1}{m+1}
\end{equation}
when $m=n$; we want to verify \eqref{momentdm} for $m=n+1$.
Since $\Gamma_{\mu}$ is Hankel, we have that $\gamma_{n+1,0}=\gamma_{n,1}$, that is,
$$
\binom{n+1}{n+1}\int_{0}^{1}t^{n+1}\,d\mu(t)=\binom{n+1}{n}\int_{0}^{1}(1-t)t^{n}\,d\mu(t)\,,
$$
which implies that
$$
    \int_{0}^{1}t^{n+1}\,d\mu(t)= (n+1)\left(\frac{1}{n+1}-\int_{0}^{1}t^{n+1}\,d\mu(t)\right)
$$
and, consequently,
$$
    \int_{0}^{1}t^{n+1}\,d\mu(t)=\frac{1}{n+2}\,.
$$
Therefore, by induction, \eqref{momentdm} holds for every $m \in \N$.
Since the Hausdorff moment problem admits a unique solution \cite[Theorem 2.6.4]{Akhiezer2021}, $\mu$ has to be the Lebesgue measure.
\end{proof}
 As a consequence of Theorem \ref{Hankel}, we  note that the operators $\G_\mu$ are distinct from the generalized Hilbert operators
\[
\cH_\mu f(z)\,=\,\int_0^1f(t)\frac{1}{1-tz}\,d\mu(t)\,,
\]
studied in the literature (see for example \cite{Chatzifountas2014, Galanopoulos2010, Girela2018}), for any measure $\mu$ other than constant multiples of the Lebesgue
measure. We refer the interested readers to the recent note \cite{Blasco2022} on generalized Hilbert operators.

Using Proposition \ref{the adjoint},  we can explicitly compute the adjoint of $\Gamma_\mu$.
\begin{lem}\label{the adjoint Gm}
Let $1<p<\infty$. The adjoint operator of $\Gamma_{\mu}$ on $H^p$, denoted by $\Gamma_{\mu}^{*}$, is given by
$$
\Gamma_{\mu}^{*}f(z)=\int_{0}^{1}T_{t}^{*}f(z)\,d\mu(t)=\int_{0}^{1}T_{1-t}f(z)\,d\mu(t)\, .
$$
\end{lem}
\begin{proof}
Let $\frac{1}{q}+\frac{1}{p}=1$. Again, we can assume that $f(z)$ and $g(z)$ are polynomials. By using Proposition \ref{the adjoint} and Fubini's Theorem, we obtain that
\begin{align*}
    \int_{\T} \Gamma_\mu f (\zeta) \overline{g(\zeta)}\,|d\z |\,&= \int_{\T} \left( \int_{0}^1T_t f (\zeta)d\mu(t)\right) \,  \overline{g(\zeta)}\,|d\z|
    \\[0.1in]
    &=\,\int_{0}^1 \int_{\T} T_t f (\zeta)\overline{g(\zeta)}\,| d\z | \, d\mu(t)
    \\[0.1in]
    &=\, \int_{0}^1 \int_{\T}f(\zeta)\overline{T_{1-t} g (\zeta)}\,| d\z | \, d\mu(t)
    \\[0.1in]
    &=\, \int_{\T}f(\zeta) \overline{ \left( \int_{0}^1 {T_{1-t} g (\zeta)} d\mu(t) \right) }\, | d\z |\,,
\end{align*}
which, since the polynomials are dense both in $H^p$ and $H^q$, implies the statement.
\end{proof}

In order to prove our next theorem, we provide an estimate of the norm of the operator $T_t$.
We first need this preliminary calculation.

\begin{lem}\label{change of var lem}
Let $1\leq p <\infty$. Then:
$$
\norm{T_{t}f}_{H^p}\,=\,\frac{t^{\frac{1}{p}-1}}{(1-t)^{\frac{1}{p}}}\left(\int_{\partial D(\frac{1}{2-t},\frac{1-t}{2-t})}|f(w)|^{p}|w|^{p-2}\, \dfrac{| dw |}{2\pi}\right)^{\frac{1}{p}}
$$
for every $f \in H^p$ and $0< t< 1$.
\end{lem}
\begin{proof}
 Let $f\in H^p$. We note that
 \begin{equation}\label{phi prime}
     \varphi_{t}(\mathbb{D})=D\left( \frac{1}{2-t},\frac{1-t}{2-t}\right) \quad \text{ and } \quad \varphi_{t}'(z)=\frac{1-t}{t}\cdot(\varphi_{t}(z))^{2}
 \end{equation}
for every $0< t< 1$ and $z\in \mathbb{D}$, where $D(z_0, r)$ is the open disc centered at $z_0$ of radius $r$. Applying a change of variables, we get
\begin{align*}
||T_{t}f||_{H^p}&=\frac{1}{t}\left(\int_{\mathbb{T}}|f(\varphi_{t}(\zeta))|^{p}|\varphi_{t}(\zeta)|^{p}\, \dfrac{| d\z |}{2\pi} \right)^{\frac{1}{p}}\\
&=\frac{t^{\frac{1}{p}-1}}{(1-t)^{\frac{1}{p}}}\left(\int_{\mathbb{T}}|f(\varphi_{t}(\zeta))|^{p}|\varphi_{t}(\zeta)|^{p-2}
|\varphi_{t}'(\zeta)|\, \dfrac{| d\z |}{2\pi} \right)^{\frac{1}{p}}\\
&=\frac{t^{\frac{1}{p}-1}}{(1-t)^{\frac{1}{p}}}\left(\int_{\partial D(\frac{1}{2-t},\frac{1-t}{2-t})}|f(w)|^{p}|w|^{p-2}\, \dfrac{| dw |}{2\pi} \right)^{\frac{1}{p}}\,.
\end{align*}
\end{proof}
We are now ready for an estimate of the norm of $T_t$ when $1\leq p<\infty.$
\begin{lem}\label{normT_t p}
Let $\,0<t<1$ and $1< p<\infty$. For the norm of $T_t$ on $H^p$ we have
$$
||T_{t}||_{H^p\to H^p}\,\leq\, B_p\,\frac{t^{\frac{1}{p}-1}}{(1-t)^{\frac{1}{p}}}
$$
for some constant $B_p>0$, depending only on $p$, which can be chosen equal to $1$ when $2\leq p<\infty$. Moreover, when $p=1$,
$$
||T_{t}||_{H^1 \to H^1}\leq C \log\left( \frac{e}{t}\right) \frac{1}{1-t}
$$
for an absolute constant $C>0$.
\end{lem}
\begin{proof}
Let us first consider the case $p\geq 2$. Then, by applying Lemma \ref{change of var lem} and \cite[Theorem 2.1]{Gabriel1928}, we get
\begin{align*}
\norm{T_t f}_{H^p}&=\frac{t^{\frac{1}{p}-1}}{(1-t)^{\frac{1}{p}}}\left( \int_{\partial D(\frac{1}{2-t},\,\frac{1-t}{2-t})}|f(w)|^{p} |w|^{p-2}\,\dfrac{| dw |}{2\pi} \right)^{\frac{1}{p}}\\
&\leq \frac{t^{\frac{1}{p}-1}}{(1-t)^{\frac{1}{p}}}\left( \int_{\bT}|f(w)|^{p}\, \dfrac{| dw |}{2\pi} \right)^{\frac{1}{p}}
\\
& = \frac{t^{\frac{1}{p}-1}}{(1-t)^{\frac{1}{p}}}\, \norm{f}_{H^p}\,.
\end{align*}

Let us now take in consideration $1< p<2$. We note that every function $f \in H^p$ can be written as $f(z)=f(0)+zg(z)$ for some $g\in H^p$.
Therefore
$$||T_t f||_{H^p}\leq |f(0)|\cdot ||T_{t}(1)||_{H^p}+||T_{t}(Sg)||_{H^p}$$
where $Sg(z)=zg(z)$ is the shift operator acting on $H^p$. For the first term we know that
\[
|f(0)|\leq 2^{\frac{1}{p}}||f||_{H^p}
\]
and
\[
\begin{split}
||T_{t}(1)||_{H^p}&=\left(\int_{-\pi}^{\pi}\frac{1}{|1-(1-t)e^{i\theta}|^{p}}\,d\theta\right)^{\frac{1}{p}}\\
&\leq\,C \left(\int_{0}^{\pi}\frac{1}{(t+\theta)^{p}}\,d\theta\right)^{\frac{1}{p}}
<\, C\,\left(\frac{1}{p-1}\right)^{\frac{1}{p}}t^{\frac{1}{p}-1}\,,
\end{split}
\]
where we have used a classical estimate (see for example \cite[Proposition 1.23]{Pavlovic2019}) for the first inequality.
For the second term, considering \eqref{phi prime} and using  \cite[Theorem 2.1]{Gabriel1928}, we obtain that
\begin{align*}
||T_{t}(Sg)||_{H^p}& = \frac{1}{t}\left(\int_{\mathbb{T}}|g(\varphi_{t}(\zeta))|^{p}|\varphi_{t}(\zeta)|^{2p}\,\dfrac{| d\z |}{2\pi} \right)^{\frac{1}{p}}\\
&=\frac{t^{\frac{1}{p}-1}}{(1-t)^{\frac{1}{p}}}\left(\int_{\mathbb{T}}|g(\varphi_{t}(\zeta))|^{p}|\varphi_{t}(\zeta)|^{2p-2}|\varphi_{t}'(\zeta)|\,\dfrac{| d\z |}{2\pi} \right)^{\frac{1}{p}}\\
&=\frac{t^{\frac{1}{p}-1}}{(1-t)^{\frac{1}{p}}}\left(\int_{\varphi_{t}\mathbb{(T)}}|g(w)|^{p}|w|^{2p-2}\,\dfrac{| dw |}{2\pi}\right)^{\frac{1}{p}}\\
&\leq \frac{t^{\frac{1}{p}-1}}{(1-t)^{\frac{1}{p}}}\;\norm{g}_{H^p}\,.
\end{align*}
This means that
\begin{align*}
||T_{t}(Sg)||_{H^p}& \leq \frac{t^{\frac{1}{p}-1}}{(1-t)^{\frac{1}{p}}}\,||Sg||_{H^p}
\\
&= \frac{t^{\frac{1}{p}-1}}{(1-t)^{\frac{1}{p}}}\,||f-f(0)||_{H^p}\\
&\leq \frac{t^{\frac{1}{p}-1}}{(1-t)^{\frac{1}{p}}}(2^{\frac{1}{p}}+1)\,||f||_{H^p}\,.
\end{align*}

The above computations imply that
$$
||T_{t}||_{H^p}\leq C \left(\frac{1}{p-1}\right)^{\frac{1}{p}}\frac{t^{\frac{1}{p}-1}}{(1-t)^{\frac{1}{p}}}\,.
$$
The proof of the case $p=1$ requires similar arguments, with very few natural modifications, so we omit the details.
\end{proof}

We are now ready to discuss our central result, concerning the boundedness of the operators $\G_\mu$.
\begin{proof}[Proof of Theorem \ref{main theom p>1} ]
First of all let us verify that \eqref{condition p>1} is necessary for the boundedness of $\Gamma_\mu$. Let $0< a < \frac{1}{p}$, so that $f_{a}(z)=\frac{1}{(1-z)^{a}} \in H^p$. Consequently
\begin{equation}\label{Gfa}
\begin{split}
    \Gamma_{\mu}(f_{a})(z)&= \bigintsss_0^1\frac{1}{\left[ 1-\frac{t}{1-(1-t)z}\right]^a}\ \frac{1}{1-(1-t)z} d\mu(t)
    \\[0.1in]
    &=\int_0^1 \frac{[1-(1-t)z]^{a-1}}{\left[(1-t)(1-z)\right]^a}d\mu(t)\\[0.1in]
    &= f_{a}(z)\,\int_{0}^{1}\frac{[1-(1-t)z]^{a-1}}{(1-t)^{a}}\,d\mu(t) \,,
    \end{split}
\end{equation}
and we note that
\[
\begin{split}
\left|\int_{0}^{1}\frac{(1-(1-t)z)^{a-1}}{(1-t)^{a}}\,d\mu(t)\right| &\geq \mbox{Re}\int_{0}^{1}\frac{(1-(1-t)z)^{a-1}}{(1-t)^{a}}\,d\mu(t)\\
&\geq \int_{0}^{1}\frac{(2-t)^{a-1}}{(1-t)^{a}}\,d\mu(t)
\\
&\geq \frac{1}{2^{1-a}}\int_{0}^{1}\frac{1}{(1-t)^{a}}\,d\mu(t)\,.
\end{split}
\]

Since $\Gamma_\mu$ is bounded in $H^p$, the above computation implies that
\begin{equation*}
\frac{1}{2^{1-a}}\int_{0}^{1}\frac{1}{(1-t)^{a}}\,d\mu(t) \, ||f_a||_{H^p} \leq ||\Gamma_{\mu}f_{a}||_{H^p}\leq ||f_{a}||_{H^p}||\Gamma_{\mu}||_{H^p\to H^p} ,
\end{equation*}
and letting $a\to \frac{1}{p}$
\begin{equation}\label{normG1}
    \frac{1}{2^{1-\frac{1}{p}}} \int_{0}^{1}\frac{1}{(1-t)^{\frac{1}{p}}}\,d\mu(t) \leq ||\Gamma_{\mu}||_{H^p\to H^p}.
\end{equation}

With similar reasoning, since $\Gamma^*_\mu$ is bounded in $H^q$ when $\frac{1}{q}+\frac{1}{p}=1$, we observe that if $0< b < \frac{1}{q}$, then by Lemma \ref{the adjoint Gm}
\begin{equation}\label{G*fa}
    \Gamma_{\mu}^{*}(f_{b})(z)\,=\,f_{b}(z)\,\int_{0}^{1}\frac{(1-tz)^{b-1}}{t^{b}}\,d\mu(t) \,,
\end{equation}
and
\begin{align*}
\left|\int_{0}^{1}\frac{(1-tz)^{b-1}}{t^{b}}\,d\mu(t)\right|\,\geq \, \frac{1}{2^{1-b}}\int_{0}^{1}\frac{1}{t^{b}}\,d\mu(t)\,.
\end{align*}
Therefore, letting $b \to \frac{1}{q}$,
\begin{equation}\label{normG2}
\frac{1}{2^{\frac{1}{p}}}\int_{0}^{1}\frac{1}{t^{1-\frac{1}{p}}}\,d\mu(t)\leq ||\Gamma_{\mu}^{*}||_{H^q\to H^q} \leq\,C_p\, ||\Gamma_{\mu}||_{H^p\to H^p}\,.
\end{equation}
Finally, putting together \eqref{normG1} and \eqref{normG2}, we have
\[\small
\begin{split}
\int_{0}^{1}\frac{t^{\frac{1}{p}-1}}{(1-t)^{\frac{1}{p}}}\, d\mu(t) & \leq
2^{\frac{1}{p}}\int_{0}^{1/2}\frac{1}{t^{1-\frac{1}{p}}}\, d\mu(t) +\left( \frac{1}{2}\right)^{\frac{1}{p}-1} \int_{1/2}^{1}\frac{1}{(1-t)^{\frac{1}{p}}}\, d\mu(t)
\\[0.1in]
&=2^{\frac{1}{p}}\int_{0}^{1}\frac{1}{t^{1-\frac{1}{p}}}\,d\mu(t) + 2^{1-\frac{1}{p}}\int_{0}^{1}\frac{1}{(1-t)^{\frac{1}{p}}}\,d\mu(t)
\\[0.1in]
&\leq \,C_p\,||\Gamma_{\mu}||_{H^p\to H^p}\,.
\end{split}
\]

On the other hand, it is clear that the condition \eqref{condition p>1} is also sufficient. Indeed, by Lemma \ref{normT_t p} and the generalized Minkowski's inequality, we have that
$$
    ||\Gamma_{\mu} f||_{H^p}\leq \int_0^1 \n{T_t f}_{H^p}\,d\mu(t)\leq B_p \,\norm{f}_{H^p} \int_{0}^{1}\frac{t^{\frac{1}{p}-1}}{(1-t)^{\frac{1}{p}}}\,d\mu(t)\ .
$$
Moreover, from Lemma \ref{normT_t p}, when $p\geq 2$ the constant $B_p$ can be chosen equal to $1$.

The case $p=1$ is treated in a similar way. Using the exact same reasoning, one can verify that the condition \eqref{condition p>1} is sufficient for the boundedness of $\Gamma_\mu$ in $H^1$.

Conversly, let us fix $f=1$. Then we have
\begin{align*}
\G_{\m}(1)(z)&=\int_0^1\frac{1}{1-(1-t)z}\,d\m(t)\\
&=\int_0^1\sum_{n=0}^{\infty}(1-t)^nz^n\,d\m(t)\\
&=\sum_{n=0}^{\infty}\left(\int_0^1(1-t)^n\,d\m(t)\right)z^n\,.
\end{align*}
Applying Hardy's inequality we get
\[
\sum_{n=0}^{\infty}\frac{1}{n+1}\int_0^1(1-t)^n\,d\m(t)\leq \pi\n{\G_{\m}(1)}_{H^1}\,,
\]
hence
\begin{align*}
\int_0^1\frac{1}{1-t}\log\left(\frac{e}{t}\right)\,d\m(t)&=
\int_0^1\frac{1}{1-t}\log\frac{e}{1-(1-t)}d\m(t)\\[0.1in]
&= \int_0^1\sum_{n=0}^{\infty}\frac{1}{n+1}(1-t)^n\,d\m(t)\,\,\\[0.1in]
&\leq\,\pi\n{\G_{\m}(1)}_{H^1}<\infty \,,
\end{align*}
which proves our claim.

Next, we obtain the exact value of the norm of the operator
\[
\G_\mu:\,H^p \to H^p,\;\; 2\leq p < \infty\,.
\]
For that range of $p$, as noted before, Lemma \ref{normT_t p} combined with Minkowski's inequality gives a sharp upper bound
\[
\|\Gamma_\mu\|_{H^p\to H^p}\leq \int_0^1 \frac{t^{\frac{1}{p}-1}}{(1-t)^{\frac{1}{p}}}\,d\mu(t)\, . 
\]
For this reason, we need to estimate a sharp, lower bound for the norm of $\Gamma_\mu$. Some of the ideas that we employ here, implicitly appeared in \cite{LINDSTROM2022} in connection to the exact value of the {\it essential norm} of a class of generalized integral operators.

For $\delta >0$, define
\[
    \Gamma_\mu^\delta(f)(z):=\int_{\delta}^{1}T_t(f)(z)d\mu(t).
\]
Let $f_{a}(z)=\dfrac{1}{(1-z)^{a}}$, $0<a<\frac{1}{p}$. We note that 
\begin{equation}\label{E_2}
    \lim_{a \to \frac{1}{p}}\|f_a\|_{H^{p}}=\infty\,,
\end{equation}
and the sequence of the normalized functions $\{\frac{f_a}{\norm{f_a}}\}_a$ tends weakly to zero as $a\to \frac{1}{p}$. For the ease of notation, we write $\Tilde{f_a}=\frac{f_a}{\norm{f_a}}$. 

Working as before, we write  $\Gamma^\delta_{\mu}(f_a)=\Lambda_{a}\cdot f_{a}$, where
$$
\Lambda_{a}(z)=\int_{\delta}^{1}\frac{(1-(1-t)z)^{a-1}}{(1-t)^{a}}\,d\mu(t).
$$

 Let $(a_{n},z_{n})\to (\frac{1}{p},1)$, where $ 0 < a_{n} < \frac{1}{p}$ for each $n\in \mathbb{N}$, and $\{a_n\}_n$ is an increasing sequence. Since, for $t \in [\delta,1]$,
 \[
 \Big|\frac{(1-(1-t)z)^{a_n-1}}{(1-t)^{a_n}}\Big|\leq \frac{\delta^{a_1-1}}{|1-t|^{\frac{1}{p}}}\, , 
 \]
by an application of the dominated convergence theorem we get
\begin{align*}
\lim_{(a_n,z_n)\to (\frac{1}{p},1)}&|\Lambda_{a_n}(z_n)|=\Lambda_{\frac{1}{p}}(1)=\int_{\delta}^{1}\frac{t^{\frac{1}{p}-1}}{(1-t)^{\frac{1}{p}}}\,d\mu(t)
\end{align*}
and, therefore, for every $z \in \overline{\mathbb D}$ we have
\begin{equation}\label{E_Lambda} 
\lim_{(a,z)\to (\frac{1}{p},1)}|\Lambda_{a}(z)|=\Lambda_{\frac{1}{p}}(1).
\end{equation}

Notice that
\[
\abs{\;\norm{\G^\delta_\mu (\Tilde{f_a})}_{H^p}\,-\,\Lambda_{\frac{1}{p}}(1)\;}^p \leq 
\norm{\Tilde{f_a}\,(\Lambda_a - \Lambda_{\frac{1}{p}}(1))}^p_{H^p}\,.
\]
Let $r>0$, and set $D_r=\{z\in \C:\,\abs{z-1}\leq r\}$. Then
\[
\begin{split}
\norm{\Tilde{f_a}\,(\Lambda_a - \Lambda_{\frac{1}{p}}(1))}^p_{H^p}\,&\leq \,
\sup_{\zeta \in \partial \D \setminus D_r}\, \abs{\Tilde{f_a}(\zeta)}^p \left( \norm{\Lambda_a}_{H^p}\,+\,\Lambda_{\frac{1}{p}}(1) \right)^p\,+
\\
& \quad +\,\sup_{\zeta \in \partial \D \cap D_r} \abs{\Lambda_a(\zeta) - \Lambda_{\frac{1}{p}}(1)}^p\,\norm{\Tilde{f_a}}^p_{H^p}\,.
\end{split}
\]
For the first term, for $a$ sufficiently close to $\frac{1}{p}$, we have that 
$$\sup_{\zeta\in \partial\D\setminus D_r}\,\abs{\Tilde{f_a}(\zeta)}^p < \epsilon\,,$$
 since $\Tilde{f_a} \to 0$, uniformly in $\abs{\zeta -1}>r$, as $a\to \frac{1}{p}$. For the second term, considering \eqref{E_Lambda}, we have that $\abs{\Lambda_a(\zeta) - \Lambda_{\frac{1}{p}}(1)}< \epsilon$, hence $\norm{\G^\delta_\mu (\Tilde{f_a})}_{H^p} \to \Lambda_{\frac{1}{p}}(1)$, as $a\to \frac{1}{p}$, which, in turn, implies that
 \[
 \norm{\G^\delta_\mu}_{H^p \to H^p} \geq \Lambda_{\frac{1}{p}}(1)\,.
 \]
 Consequently, 
 \begin{align*}
\|\Gamma_\mu\|_{H^p\to H^p}\,=&\;\|\Gamma_\mu-\Gamma_{\mu}^\delta+\Gamma_{\mu}^\delta\|_{H^p\to H^p}\\
\geq& \;  \|\Gamma_{\mu}^\delta\|_{H^p\to H^p} - \|\Gamma_\mu-\Gamma_{\mu}^\delta\|_{H^p\to H^p}\\
\geq&\; \int_\delta^1 \dfrac{t^{\frac{1}{p}-1}}{(1-t)^{\frac{1}{p}}}\,d\mu(t) -\int_0^\delta \dfrac{t^{\frac{1}{p}-1}}{(1-t)^{\frac{1}{p}}}\,d\mu(t) ,
\end{align*}
and letting $\delta \to 0^+$, the result follows.
\end{proof}

As a final result for this section, we study the compactness of the operators $\G_\mu$.

\begin{proof}[Proof of Theorem \ref{Compact}]
Let $1<p<\infty$ and $\frac{1}{q}+\frac{1}{p}=1$. Suppose, on the contrary, that $\Gamma_{\mu}$ is compact on  $H^p$. Since the family of functions
\[
 k_{w}(z)=\frac{(1-|w|^{2})^{\frac{1}{q}}}{1-\overline{w}z}
\]
converges to zero weakly as $|w|\to 1$, we have that $\norm{\Gamma_{\mu}(k_w)}_{H^p} \to 0$.
By considering the classical growth estimate \eqref{growth} for functions in $H^p$ , we get
\begin{align*}
\left|\int_{0}^{1}\frac{(1-|w|^{2})^{\frac{1}{q}}(1-|z|^{2})^{\frac{1}{p}}}{1-(1-t)z-\overline{w}t}\,d\mu(t)\right|&=
(1-|z|^{2})^{\frac{1}{p}}\left|\Gamma_{\mu}(k_w)(z)\right|\\
&\leq 2^{\frac{1}{p}}||\Gamma_{\mu}(k_w)||_{H^p}\,.
\end{align*}
Let $z=w=r \in (0,1)$. Then
\begin{align*}
\mu(0,1)\leq\frac{1-r^{2}}{1-r}\mu(0,1)
\leq 2^{\frac{1}{p}}||\Gamma_{\mu}(k_r)||_{H^p}\,.
\end{align*}
By letting  $r\to 1^{-}$, we obtain that $\mu(0,1)=0$, which is a contradiction. Consequently $\Gamma_\mu$ cannot be compact in $H^p$,  $p>1$.

Next, we deal with the case $p=1$. Suppose that for some measure $\mu$ , $\Gamma_\mu$ is compact on $H^1$.  Let $k_r(z)=\frac{1}{1-rz},\, 0<r<1$ and $z\in \mathbb{D}.$ Notice that $\Vert k_r \Vert_{H^1}$ is comparable to $ \log\big(\frac{e}{1-r}\big)$, as $r\to 1^{-}$, and consider the normalized functions $ \hat{k}_r = k_r / \Vert k_r \Vert_{H^1}$. Then the compactness of $\Gamma_\mu$ implies that $\Gamma_\mu(\hat{k}_r) $ converges strongly to zero \cite[Lemma 3.7]{Tjani2003} as $r\to 1^-$. We shall use the Fej\'er-Riesz inequality in order to arrive at a contradiction. Notice that on the segment $(0,1)$ the functions $\Gamma_\mu(\hat{k}_r)$ are positive, therefore for some positive constant $C$, we get

  \begin{align*}
      \int_0^1|\Gamma_\mu(\hat{k}_r)(s)| ds & = \frac{1}{\Vert k_r \Vert_{H^1}}  \int_0^1 \int_0^1 \frac{1}{1-rt-(1-t)s} d\mu(t)  ds \\
      = & \frac{1}{\Vert {k}_r \Vert_{H^1}} \int_0^1 \frac{1}{1-t} \log\bigg( \frac{1-rt}{t(1-r)} \bigg) d\mu(t) \\
      & \geq\, C \int_0^1\frac{1}{1-t} \log\bigg( \frac{1-rt}{t(1-r)} \bigg)\log\bigg( \frac{e}{1-r} \bigg)^{-1} d\mu(t)\,.
  \end{align*}

  Now applying Fatou's Lemma, as $r\to 1^-$, we find that
  \[ 
  C \int_{0}^1 \frac{1}{1-t}d\mu(t) \leq   \liminf_{r\to 1^- } \int_0^1 |\Gamma_\mu(\hat{k}_r)(s)| ds \,\leq\, \pi   \lim_{r\to 1^-} \Vert \Gamma_\mu(\hat{k}_r) \Vert_{H^1} = 0\,. 
  \]
  This implies that $\mu = 0$.

 To show that $\Gamma_\mu$ is completely continuous whenever it is bounded, notice that $T_t$ is completely continuous for every $t\in (0,1)$, by a theorem of Cima and Matheson \cite[Proposition 1]{Cima1994}. Consider now a sequence of functions $\{ f_n \}$ which is weakly null in $H^1$. Then, by the complete continuity of $T_t$ we have that $\lim_n \Vert T_t f_n \Vert_{H^1} = 0$, for all $0<t<1 $. Furthermore, using Lemma \ref{normT_t p} and for some constant $C>0$, we have
  \begin{equation*}
      \Vert T_t  f_n \Vert_{H^1} \leq \sup_{n} \Vert f_n \Vert_{H^1} \Vert  T_t \Vert_{H^1} \,\leq \,C\, \sup_n \Vert f_n \Vert_{H^1} \log \frac{e}{t}\frac{1}{1-t}\,.
  \end{equation*}
  By hypothesis this is a function in $L^1((0,1), \mu)$ therefore applying the dominated convergence theorem we conclude that 
  \[
  \limsup_n \Vert \Gamma_\mu f_n \Vert_{H^1} \leq \limsup_n \int_0^1 \Vert T_t f_n \Vert_{H^1}\, d\mu(t) = 0\,.
  \]
\end{proof}



\section{Concluding remarks}
 There are still some questions regarding $\G_\mu$ which we have not been able to resolve. In specific, one might be tempted to think that the norm of $\G_\mu$ acting on $H^p,\;1<p<2$, is equal to $ \int_{0}^{1}\frac{t^{\frac{1}{p}-1}}{(1-t)^{\frac{1}{p}}}\,d\mu(t)$, as is the case for $p\geq 2$. Unfortunately, if the measure $\mu$ is a Dirac mass $\delta_{t_0} $ at $t_0\in (0,1)$, the operator $\G_\mu$ coincides with the weighted composition operator $T_{t_0}$. It is then easily verified that this is not even asymptotically true as $t_0\to 0^+$ and $p\to 1^+$.

Regarding the range $0<p<1$, one can adapt our methods to show that a necessary condition for the boundedness of the operator $\G_\mu$ on $H^p$ is that 
\[
\int_0^1 \frac{1}{(1-t)^{\frac{1}{p}}} d\mu(t) < + \infty\,.
\]
However, it is not clear whether the above condition is also sufficient for the bounedness of $\G_\mu$.

\subsection{Acknowledgments}
The authors would like to express their gratitude to professors Petros Galanopoulos and Aristomenis Siskakis, for introducing us to this problem and for all their valuable help during the preparation of this article.
We would also like to thank professors Santeri Miihkinen and Jani Virtanen for the valuable discussions on the topic, during their visit in Thessaloniki.
\bibliographystyle{plain}
\bibliography{Literature}

\end{document}